\documentclass{amsart}
\usepackage{graphicx}
\usepackage{amssymb}
\usepackage{amsthm}
\usepackage{mathrsfs,txfonts}
\usepackage[all]{xy}
\usepackage{enumerate}
\usepackage{color}
\DeclareMathAlphabet{\mathpzc}{OT1}{pzc}{m}{it}

\def\spec{\operatorname{Spec}}
\newcommand{\Oh}{\mathcal{O}}
\newcommand{\pf}{\noindent{\bf Proof:\ \ }}
\newcommand{\cqd}{{\hfill $\rule{2mm}{2mm}$}\vspace{3mm}}

\def\mult{\operatorname{mult}}
\def\min{\operatorname{min}}

\def\deg{\operatorname{deg}}
\def\car{\operatorname{char}}

\def\lim{\operatorname{lim}}
\def\I{\operatorname{I}}

\def\C{\mathbb{C}}

\def\Ob{{\overline{\mathcal{O}}}}

\def\M{{\mathfrak{m}}}

\def\mult{{\rm mult}}

\providecommand{\deg}{\mathop{\rm deg}\nolimits}

\providecommand{\ord}{\mathop{\rm ord}\nolimits}

\newtheorem{theorem}{Theorem}[section]
\newtheorem{lemma}[theorem]{Lemma}
\newtheorem{corollary}[theorem]{Corollary}

\newtheorem{proposition}[theorem]{Proposition}
\newtheorem{example}[theorem]{Example}
\newtheorem{remark}[theorem]{Remark}

\title{The Milnor number of plane branches with tame semigroups of values}
 \author{A. Hefez, \ J.H.O. Rodrigues \ and \ R. Salom\~ao \\ \\ \bigskip
{\tiny Universidade Federal Fluminense - Niter\'oi} }
\begin{document}

\begin{abstract}
The Milnor number of an isolated hypersurface singularity, defined as the
codimension $\mu(f)$ of the ideal generated by the partial derivatives of a
power series $f$ that represents locally the hypersurface, is an important
topological invariant of the singularity over the complex numbers. However it
may loose its significance when the base field is arbitrary. It turns out
that if the ground field is of positive characteristic, this number depends upon
the equation $f$ representing the hypersurface, hence it is not an invariant of
the hypersurface.
For a plane branch represented by an irreducible convergent power series $f$ in two indeterminates over the complex numbers,
it was shown by Milnor that $\mu(f)$ always coincides with the conductor $c(f)$ of the semigroup of
values $S(f)$ of the branch. This is not true anymore if
the characteristic of the ground field is positive. In this paper we show that, over algebraically closed fields
of arbitrary characteristic, this is true, provided that the semigroup $S(f)$ is tame, that is,
the characteristic of the field does not divide any of
its minimal generators.  \medskip

{\small \noindent Keywords: {Singularities in positive characteristic, Milnor number in positive characteristic, Singularities of algebroid curves}\medskip

\noindent 2010 Mathematics Subject Classification: {14B05, 14H20, 14D06, 32S05}}

\end{abstract}
\maketitle

\section{Introduction}
\label{I}

Several aspects of the local study of singularities of algebraic varieties may be reduced to the study of {\em algebroid varieties}. In particular, by a singular algebroid plane curve we mean a scheme $C=\spec(\Oh)$, where $(\Oh,\M)$ is a local one-dimensional complete $k$-algebra, with $k$ any algebraically closed field, such that $\dim\M/\M^2=2$. From the completeness of $\Oh$, for any choice of generators $x,y$ of the ideal $\M$, there is a surjection $\varphi\colon k[[X,Y]] \to \Oh$, where $k[[X,Y]]$ is the ring of formal power series in two indeterminates with coefficients in $k$. The kernel of this surjection is a principal ideal $\langle f\rangle$, which generator $f$ is uniquely determined up to a multiplication by a unit in $k[[X,Y]]$, which we call an equation of $C$. If we define $\Oh_f=k[[X,Y]]/\langle f\rangle$, and write $C_f=\spec(\Oh_f)$, we have that $\Oh\simeq \Oh_f$ and $C_f$ is isomorphic to $C$ as a scheme over $k$. When $C$ is an integral scheme, we call it a {\em plane branch}.

Two algebroid curves $C_1=\spec(\Oh_1)$ and $C_2=\spec(\Oh_2)$ will be considered equivalent, writing $C_1\sim C_2$, when they are isomorphic as $k$-schemes, that is, when $\Oh_1\simeq \Oh_2$ as $k$-algebras. In this case, if $f_1$ is an equation of $C_1$ and $f_2$ is an equation of $C_2$, then $f_1$ and $f_2$ are related by the existence of an automorphism $\Phi$ and a unit $u$ of $k[[X,Y]]$ such that $f_2=u(f_1\circ \Phi)$. In this situation, we say that $f_1$ and $f_2$ are {\em contact equivalent}.

Given $f\in k[[X,Y]]$, the ideal $T(f)=\langle f, f_{X}, f_{Y} \rangle$ is called the \emph{Tjurina ideal} of $f$ and the dimension $\tau(f)$ of the $k$-vector space $k[[X,Y]]/T(f)$ is the so called {\em Tjurina number} of $f$. It is easy to check that this number is invariant under contact equivalence, so it defines an invariant of $C$, denoted by $\tau(C)$ or $\tau(\Oh)$.

The ideal $J(f)=\langle f_X, f_Y \rangle$ is called  the {\em Jacobian ideal} of $f$ and plays an important role when $k=\C$. The {\em Milnor number} $\mu(f)$ of $f$ is the dimension as $k$-vector space of $k[[X,Y]]/J(f)$. From the chain rule it is immediate to verify that $\mu(f)=\mu(\Phi(f))$ for any automorphism $\Phi$ of $k[[X,Y]]$.

When $f\in \C\{X, Y\}$, the ring of convergent complex power series, Milnor proves by
topological methods that $\mu(uf)=\mu(f)$ for any unit $u$ in $\C\{X,Y\}$.  In arbitrary
characteristic this does not hold as one can see from the very simple example below.

\begin{example}\label{firstexample} Let $\car\,{k}=p$ and $f=Y^p+X^{p+1}\in k[[X,Y]]$. Then $C_f$ is a plane branch such that $\tau(C_f)=p^2$ and $$
\mu(f)=\infty,  \ \ \quad  {\rm while} \ \ \quad \mu((1+Y)f)=p^2 \neq \mu(f).$$
\end{example}

From the inclusion $J(f) \subset T(f)$, it is clear that $\tau(f) \leqslant \mu(f)$. So, one always has
$\mu(f) < \infty \ \Rightarrow \tau(f) < \infty$.
In characteristic zero, one also has the converse of the above implication. In positive characteristic, this converse may fail, as one can see from Example \ref{firstexample} above. \medskip

Since the ideal that defines the singular locus of a curve $C=C_f$ is its Tjurina ideal $T(f)$, it is natural to say that a curve $C$ has an {\em isolated singularity} if $0<\tau(C)<\infty$.

Notice that this is a well posed definition, since $\tau(f)=\tau(g)$ when $f$ and $g$
are contact equivalent. So, in characteristic zero, to say that $C_f$ has an isolated
singularity is equivalent to say that $0<\mu(f)<\infty$, but not in arbitrary
characteristic.

There is an easy criterion in arbitrary characteristic (cf. \cite[Proposition 1.2.11]{B})
for a plane curve $C_f$ to have an isolated singularity:

{\em $C_f$ has an isolated singularity if and only if $f$ is reduced.}

 In contrast, the fact that
$f$ is reduced is not sufficient to guarantee that $\mu(f)<\infty$ as shows Example \ref{firstexample}.
Also, the vanishing of one of the partial derivatives of $f$ implies $\mu(f)=\infty$,
but this is not a necessary condition, as the
following example shows.

\begin{example}\label{tresretas} Let $\car\,{k}=3$ and $f=X^2Y+Y^2X\in k[[X,Y]]$. We have that $f=XY(X+Y)$ is reduced, but $f_X$ and $f_Y$ are both nonzero and have the common factor $Y-X$, implying that $\mu(f)=\infty$.
\end{example}

The preceding examples highlight the need of a better notion of Milnor number and
this was treated in our previous work \cite[Section 3]{HRS} in the more general context of hypersurfaces. There we were led to define the {\it Milnor number} of $C= C_f$ as $\mu(C)=e_0(T(f))$, where $e_0(-)$ stands for the Hilbert-Samuel
multiplicity of an ideal. This number is invariant in the contact class of $f$. Also, a criterion was given there to recognize if an isolated hypersurface singularity $C$ satisfies $\mu(C)=\mu(f)$ for any of its
equations $f$. Namely it is shown in \cite[Theorem 4.8]{HRS} that this holds if and
only if $f^{\ell}\in\mathcal M T(f)^{\ell}$ for some positive integer $\ell$, where $\mathcal M$ is the maximal ideal of $k[[X_1,\ldots,X_n]]$.

Our aim in the present work is to study plane branches over arbitrary algebraically closed fields. In characteristic zero, the Milnor number $\mu(f)$ coincides with the conductor $c(f)$ of the semigroup of values $S(f)$ of the branch $\Oh_f$ (cf. \cite{Mi}), where the semigroup of values of $\Oh$ is the set of all finite intersection multiplicities of $f$ with plane algebroid curves and $c(f)$ is the smallest integer $\alpha\in S(f)$ such that $\alpha +\mathbb N \subset S(f)$. Since these notions are invariant in the contact class of $f$, then we may define $S(C)=S(f)$ and $c(C)=c(f)$ and have that $$C_1 \sim C_2 \ \Rightarrow \ S(C_1)=S(C_2) \quad \text{and}  \quad c(C_1)=c(C_2). $$

In arbitrary characteristic, Deligne proved in \cite{De} (see also \cite{MH-W}) the inequality $\mu(f)\geqslant c(C_f)$ and shows that the difference $\mu(f)-c(C_f)$ is the number of {\em wild vanishing cycles} associated to $f$.

Our main result is the proof that Milnor's number $\mu(f)$ of $f$ and the
conductor $c(C_f)$ of a branch $C_f$ coincide when the characteristic does not divide
any of the minimal generators of the semigroup of values $S(C_f)$. If this is so, we call
$S(C_f)$ a \emph{tame semigroup}. Our proof was inspired by a result of P. Jaworski in the
work \cite{Ja2}, which we simplified and extended to arbitrary characteristic,
under the appropriate assumptions. We would like to point out that in the process
of writing our results, E. Garc\'ia-Barroso and A. Ploski published the paper
\cite{GB-P}, where they show, by other methods, a weaker version of our main result (with the converse),
in the particular case when $p$ is greater than the multiplicity $\mult(f)$ of $f$. They also observed that (contrary to ours) their proof fails when $p \leq \mult(f)$. The proof we give of the result, without the restriction $p> \mult(f)$, is more subtle.

We should also mention that
H.D. Nguyen in \cite{Ng} has shown, in the case of several branches, a
generalization of this equality, but also under a strong restriction on
the characteristic. In the irreducible case, the restriction is
$p> c(C_f)+\mult(f)-1$.

This work is part of the PhD Thesis of the second author, under the supervision of the other two authors.






 \section{A fundamental formula}\label{chapIII3}

Let $\Oh$ be the ring of a singular plane branch $C$ with maximal ideal $\M$.
Consider the integral closure $\Ob$ of $\Oh$ in its field of fractions. The ideal
$\mathcal{C}(\Oh)=(\Oh:\Ob)$ is called the conductor ideal of $\Ob$ in $\Oh$. This is the largest common
ideal of $\Oh$ and $\Ob$. As an ideal of $\Ob\simeq k[[t]]$ we have
$ t^c\Ob=\mathcal{C}(\Oh)\, $ for some $c\in\mathbb{N}$. This number $c$ coincides with the
conductor $c(C)$ of the semigroup of values $S(C)$.

Choose generators $x,y$ of $\mathfrak{m}$ and a uniformizing parameter $t\in\Ob$ such that the images $x(t)$ and $y(t)$ in $\Ob\simeq k[[t]]$ are a primitive parametrization
for any equation $f$ of $C$ determined by the kernel of the epimorphism $k[[X,Y]]\rightarrow\Oh$ given by $X\mapsto x$ and $Y\mapsto y$, that is $k((x(t),y(t)))=k((t))$.
We say that $z\in\M$ is a \emph{transversal} parameter for $\Oh$ if $v(z)=\min\{v(w)\,|w\in\M \}$ where $v=\ord_t$ is the natural valuation of $\Ob$, which coincides with the intersection index with $f$, in the sense that $v(g)=\I(f,G)$, where $G\in k[[X,Y]]$ is any representative of the residual class $g\in \Oh$.
This minimum is called the multiplicity of $C$ and is denoted $\mult(C)$. We also say that $z$ is a \emph{separable} parameter with respect to $t$ if $z'(t)=\frac{dz}{dt}\neq 0$. For $h\in k[[X,Y]]$, we will write $f_x$ for $f_X(x,y)\in \Oh$ and similarly for $f_y$.

The following fundamental formula is classical and is attributed to D. Gorenstein. In order to have a proof in arbitrary
characteristic, for the convenience of the reader,  we will reproduce the one found in \cite[Appendix, page 101]{Ga}, as communicated to us by K. O. St\"ohr.

\begin{theorem} (\cite[Theorem 12]{Go})\label{Gorenstein} Let $\Oh$ be the ring of an algebroid plane branch $C$. If $x,\,y$ is any system of generators of $\M\subset\Oh$ and $f$ is an equation of $C$ corresponding to the embedding determined by $x$ and $y$, then \[v(f_y)=c(f)+v(x').\]
\end{theorem}

\pf The proof will be by induction on the number of blowing-ups necessary to dessingularize the branch $C$. If $\Oh$ is already non-singular then the result is obvious.

We treat first the case in which $x$ is both a transversal and separable parameter and the tangent line is $Y$.
Since in this coordinate system the tangent cone of an equation $f$ is $Y^n$, where
$n=\mult(C)$, then, if $f^{(1)}$ is the strict
transform of $f$, we have $X^nf^{(1)}(X,Z)=f(X,Y)$ where $Y=XZ$. Differentiation
with respect to $Z$ in the last relation leads to
\[X^nf_Z^{(1)}(X,Z)=Xf_Y(X,Y).\] This relation modulo $\langle f\rangle$,  after cancelling $X$, gives
$x^{n-1}f_z^{(1)}=f_y.$ Hence
\[v(f_y)=v(x^{n-1})+v(f_z^{(1)})=n(n-1)+v(f_z^{(1)}).\]  By the induction
hypothesis $v(f_z^{(1)})=v(x')+c(f^{(1)})$. According to a well known formula
(see for instance \cite[Formula 6.10]{He}), we have that $c(f)=c(f^{(1)})+n(n-1)$, which gives us the result in this case.

Suppose now that $x,\,y$ is an arbitrary set of generators of
$\M$. It is clear that there exists a system of generators $z,\,w$ of $\M$ such
that \[x=az+bw,\ \ \ \ y=cz+dw\] where $a,b,c,d\in k$ with $D=ad-bc\neq 0$ and,
say, $z$ is both a transversal and separable parameter of $\Oh$ and $W$ is the
tangent line. Then, for any equation $g\in k[[X,Y]]$ corresponding to the embedding
determined by $z$ and $w$, we get $g_w=f_x b+f_y d.$ It follows that
\[f_yDz'=(dx'-by')f_y=dx'f_y+bx'f_x=x'g_w,\] therefore,
$v(x')+v(g_w)=v(z')+v(f_y)$.

Since by the preceding case we have $v(g_w)=c(g)+v(z')$ and since $c(g)=c(f)$, this leads to
$v(f_y)=c(f)+v(x')$ in this general case. \cqd

We will need another formulation for Theorem \ref{Gorenstein}, which is quoted in the literature as \textit{Delgado's Formula} (cf. \cite[Proposition 7.4.1]{Ca}), proved over $\C$, which extends naturally to arbitrary algebraically closed fields.

To begin with, let $\Oh$ be the ring of an algebroid plane branch $C$ and choose an equation $f$ for $C$. For $g\in\Oh_f$, take a representative $G\in k[[X,Y]]$ of $g$ and define

$$[f,g]\colon=f_x\overline{G_y}-f_y\overline{G_x}\in\Oh\,\quad \mbox{and}\,\quad g'(t)\colon=\frac{d\,}{dt}G(x(t),y(t)).$$ Notice that $[f,g]$ and $g'(t)$ are well defined since they do not depend on the representative $G$ of $g$.   With these notations we have

\begin{corollary}\label{outrodelgado}
Let $\Oh$ be the local ring of a plane algebroid branch $C$. Fix a uniformizing parameter $t$ for $\Oh$ and let $f$ be an equation for  $C$.
Then, for every $g\in\Oh$, one has $$v([f,g])=c(f)+v(g'(t)).$$
\end{corollary}
\pf
Let $(x(t),y(t))$ be a local parametrization of $C$.
From $f(x(t),y(t))=0$ we obtain $f_X(x(t),y(t))x'+f_Y(x(t),y(t))y'=0$, hence it
is easy to check that
$$x'(f_YG_X-f_XG_Y)(x(t),y(t)) = f_Y(x(t),y(t))(g'(x(t),y(t))).$$
We deduce the result computing the orders in the preceding equality and
using Theorem \ref{Gorenstein}.\cqd

\begin{remark} Since $v(f_x)=c(f)+v(y')\geqslant c$ and $v(f_y)=c(f)+v(x')\geqslant c$,
we conclude that $$\mathrm{j}=\langle f_x,f_y\rangle\subset \mathcal{C}(\Oh).$$

Since $\Oh$ is Gorenstein, this implies that
$$\tau(\Oh) = \dim_k \Oh/j \geqslant \dim_k \Oh /\mathcal{C}(\Oh) =
\dim_k \overline{\Oh}/\Oh = \frac{c}{2}.$$
\end{remark}

\begin{corollary}\label{maisumdelgado} If  $p=\car{k}$ and $f\in k[[X,Y]]$ is irreducible, then
$$v([f,g])\geqslant c(f)+v(g)-1,$$ with equality holding if and only if $p\nmid v(g)$.
\end{corollary}

\pf This follows from the previous Corollary and from
$v(g)-1\leqslant v\Big(g'(t)\Big)$ where equality holds if and only if $p\nmid v(g)$.
\cqd




\section{Milnor number for plane branches}
\label{chapIII2}

For the definitions and notation used in this section we refer to \cite{He} where these notions are characteristic free. Let $f \in \mathcal M \subset k[[X,Y]]$ be an irreducible power series, where $k$ is an algebraically closed field of characteristic $p\geqslant 0$ and $\mathcal M=\langle X,Y\rangle$. Let us denote by $S(f)=\langle v_0,\ldots,v_g\rangle$ the semigroup of values of the branch $C_f$, represented by its minimal set of generators, where the number $g$ is called the {\em genus} of $S(f)$. These semigroups have many special properties that will be used throughout this section and which we describe briefly below.

Let us define $e_0=v_0$ and denote by $e_i=\gcd\{v_0,\ldots,v_{i}\}$ and by $n_i=e_{i-1}/e_i$, when $i=1,\ldots,g$. The semigroup $S(f)$ is {\em strongly increasing}, which means that $v_{i+1}>n_iv_i$, for $i=0,\ldots,g-1$, (cf. \cite[6.5]{He}). This implies that the
the sequence $v_0,\ldots,v_g$ is {\em nice}, which means that
$n_iv_i\in \langle v_0,\ldots,v_{i-1}\rangle$, for $i=1,\ldots,g$
(cf. \cite[Proposition 7.9]{He}). This, in turn, implies that the semigroup
$S(f)$ has a {\em conductor} $c(f)$, which is given by the formula
(cf. \cite[7.1]{He})
\[
c(f)=\sum_{i=1}^g (n_i-1)v_i-v_0+1.
\]

The semigroup $S(f)$ is also symmetric (cf. \cite[Proposition 7.7]{He}), that is,
\[
\forall\,  z\in \mathbb N, \ z\in S(f) \ \Longleftrightarrow \ c(f)-1-z \not\in S(f).
\]
Moreover, any element $x\in S(f)$ may be written in a unique way as
$x=\sum_{i=0}^{g}x_iv_i$ with
$x_{0}\in {\mathbb N}$ and $0\leqslant x_i\leqslant n_{i}-1$. We refer to this
representation as the {\em canonical representation} of $x$.

Recall that we defined a semigroup of values of a plane branch to be tame if $p=\car(k)$ does not
divide any element in its minimal set of generators.

A first property of curves with tame semigroups is the following:

\begin{lemma} \label{tameimplicafinito} If $C=\spec(\Oh)$ is a branch with tame semigroup $S(C)$, then, for any equation $f$ of $C$, one has that $\mu(f)<\infty$.
 \end{lemma}

\pf Let $f$ be an equation for $C$. We have that $\mu(f)=\infty$ if and only if there is an irreducible
$q\in k[[X,Y]]$ that divides both $f_X$ and $f_Y$, say $f_X=Aq$ and $f_Y=Bq$, for some $A,\,B\in k[[X,Y]]$. If
this is so, we write the canonical representation of the intersection multiplicity $\I(f,q)\in S(f)\setminus \{0\}$ as $\sum_{i=\ell}^g\,x_iv_i$, $x_i\in\mathbb{N}$, where
$\ell\in\{0,\ldots,g\}$ is the least index $i$ for which $x_i>0$ and we take
$h \in k[[X,Y]]$ such that
$\I(f,h)=v_{\ell}$. Then
$[f,h]=q \left(A h_Y-B h_X\right)$, and hence by computing
intersection multiplicities with $f$, and using that $S(f)$ is tame, we get that
\[v_{\ell}+c(f)-1=\sum_{i=\ell}^g\,x_iv_i \, +\gamma, \,\,
\mbox{with}\,\, \gamma\in S(f).\] Since $x_{\ell}>0$, this gives the absurd equality
(the summation is taken to be zero if $\ell=g$)
\[c(f)-1=(x_{\ell}-1)v_{\ell}+\sum_{i=\ell+1}^g\,x_iv_i \, +\gamma\in S(f).\]
 \cqd

Recall that two plane analytic branches over the complex numbers are {\em equisingular} if their semigroups of values coincide. We will keep this terminology even in the case of positive characteristic.

The following example will show that $\mu(C_f)=e_0(T(f))$ may not be constant in an
equisingularity class of plane branches.

\begin{example} The curves given by $f=Y^3-X^{11}$ and $h=Y^3-X^{11}+X^8Y$ are
equisingular with semigroup of values $S=\langle 3,11\rangle$. In
characteristic $3$, one has that $\mu(C_f)=\mu((1+Y)f)=30$ and
$\mu(C_h)=\mu((1+X)h)=24$. Notice that in this
case $S$ is not tame.
\end{example}

The following is an example which shows that the $\mu$-stability is not a
character of an equisingularity class.

\begin{example} Let $S=\langle 4,6,25\rangle$ be a strongly increasing semigroup
with conductor $c=28$. Consider the equisingularity class determined by $S$ over
a field of characteristic $p=5$. If $f=(Y^2-X^3)^2-YX^{11}$, which belongs to
this equisinsingularity class, we have that $\mu(f)=41$ and $\mu(C_f)=30$,
hence $C_f$ is not $\mu$-stable. But from \cite[Theorem 4.8]{HRS}, the
equisingular curve with equation
$h=(Y^2-X^3+X^2Y)^2-YX^{11}$ is $\mu$-stable, since $h^3\in \mathcal M T(h)^3$.
In this case one has $\mu(h)=\mu(C_h)=29$. Notice that here, again, $S$ is
not tame.
\end{example}

We now state our main result, which proof will occupy the rest of the paper.

\begin{theorem}[Main Theorem]\label{tame} If $C=\spec(\Oh)$ is a plane singular branch with $S(C)$ tame, then
$\mu(f)=\mu(C)=c(f)$ for any equation $f$ of $C$. In particular,
$C$ is $\mu$-stable.
\end{theorem}

The proof we give of this theorem is based on the following theorem which was
stated without a proof over the complex numbers in \cite{Ja1}, but proved in the
unpublished work \cite{Ja2}. Our proof, in arbitrary characteristic, is inspired
by that work, which we suitably modified in order to make it work in the more
general context we are considering.

\begin{theorem}[Key Theorem] \label{theorem1}
Let $f\in \mathcal M^2\subset k[[X,Y]]$ be an irreducible Weierstrass polynomial such that $S(f)$ is tame.
Then any family $\mathcal{F}$ of elements inside $k[[X]][Y]$ of degree in $Y$
less than $\mult(f)$ such that
\[
\{\I(f,h); \ h\in \mathcal{F}\}= S(f)\setminus \ (S(f)+c(f)-1)
\]
is a representative set of generators of the $k$-vector space $k[[X,Y]]/J(f)$.
\end{theorem}

We postpone the proof of this theorem until the next section, since it is long
and quite technical.

In order to apply Theorem \ref{theorem1} to prove Theorem \ref{tame} we need a
process that transforms a power series into a Weierstrass polynomial. This is
classically done by using Weierstrass Preparation Theorem, but this is not
appropriate for studying Milnor's
number in positive characteristic, since it involves the multiplication of the
power series by a unit and this affects the Milnor number. So, we will need a
preparation theorem that involves only coordinate changes and this will be done using a
result due to N. Levinson (cf. \cite{Le}) over $\C$, which we state so that the
same proof remains valid over arbitrary algebraically closed fields.
\medskip

\noindent {\bf Levinson's Preparation Theorem}
{\it Let
$0\neq f(X,Y)\in\mathcal M\subset k[[X,Y]]=k[[X,Y]]$,
where $k$ is an algebraically closed field of characteristic $p\geqslant 0$.
Write
\[f=\sum_{i,j}\alpha_{i,j}X^iY^j\] where $(i,j)\in\mathbb{N}^2$,
and suppose that $f$ contains for some integer $r>1$ a monomial $Y^r$ with
nonzero coefficient. If $r$ is minimal with this property and $p$ does not divide
$r$, then there exists a change of coordinates $\varphi$ in $k[[X,Y]]$, preserving
$k[[X]]$, which transforms $f$ into
\[\varphi(f)=A_0(X)\,Y^r+A_{1}(X)\,Y^{r-1}+\cdots
+A_{r-1}(X)\,Y+A_r(X),\] where
$A_i(X)\in k[[X]]$ for every $i$ and
\[A_{1}(0)=\cdots=A_{r-1}(0)=A_r(0)=0, \ \ \ A_0(0)\neq 0.\]}

For the moment we observe here the following refinement of the above
Theorem for the case of plane branches.

\begin{corollary} Let $f\in k[[X,Y]]$ be irreducible where $k$ is algebraically
closed of characteristic $p$. Set $n=\mult(f)$. If $p\nmid n$, then there exists an
automorphism $\varphi$ of $k[[X,Y]]$ such that
\[\varphi(f)=Y^n+B_{1}(X)Y^{n-1}+\cdots+ B_{n-1}(X)Y+B_n(X),\]
where
$B_i(X)\in k[[X]]$ and $\mult(B_{i}(X))>i$, for all $i=1,\ldots,n$.
\end{corollary}
\pf Since $f$ is irreducible, we have that $f=L^n+hot$, where $L$ is a linear
form in $X$ and $Y$. By changing coordinates, we may assume that $f$ is as in
the conclusion of Levinson's Preparation Theorem. Now, since $p\nmid n$, we take
an $n$-th root of $A_0(X)$ and perform the change of coordinates
$Y\mapsto YA_0^{\frac{1}{n}}$ and $X\mapsto X$. So, after only changes of
coordinates $\varphi$, we have that
\[\varphi(f)=Y^n+B_{1}(X)Y^{n-1}+\cdots+ B_{n-1}(X)Y+B_n(X),\]
is a Weierstrass polynomial, that is,
$\mult(B_{i}(X)) > i$, for $i=1,\ldots, n$.
\cqd

\noindent{\bf Proof of Theorem \ref{tame}:} From Deligne's results in \cite{De}
(see also \cite{MH-W}) and from Lemma \ref{tameimplicafinito} one always has
$c(f) \leqslant \mu(f) <\infty$.

Now, after a change of coordinates, that does not affect the result, we may
assume that the equation $f$ of $C$ is a Weierstrass polynomial. For every
$\alpha \in S(f)\setminus (S(f)+c(f)-1)$, take an element $g\in k[[X,Y]]$ such
that $\I(f,g)=\alpha$ and after dividing it by $f$ by means of the Weierstrass
Division Theorem, we get in this way a family $\mathcal F$ as in
Theorem \ref{theorem1}.

Theorem \ref{theorem1} asserts that the residue classes of the elements in
$\mathcal{F}$ generate $k[[X,Y]]/J(f)$, hence
$\mu(f)\leqslant \# \left(S(f)\setminus (S(f)+c(f)-1)\right)$.
The result will then follow from the next Lemma that asserts that the number in
the right hand side of the inequality is just $c(f)$.

The $\mu$-stability follows from the fact that for every invertible element $u$
in $k[[X,Y]]$,
both power series $f$ and $uf$ can be individually prepared into
Weierstrass form by means of a change of coordinates that does not alter the
semigroup, nor the Milnor numbers. Hence, $\mu(f)=c(f)=c(uf)=\mu(uf).$ \cqd

\begin{lemma} \label{contagem}
 $\# \left(S(f)\setminus (S(f)+c(f)-1)\right)=c(f)$.
\end{lemma}
\pf In fact, to every $i\in\{0,1,\ldots,c(f)-1\}$ we associate
$s_i\in S(\Oh)\setminus (S(\Oh)+c(\Oh)-1)$ in the following way:
\[
s_i=\left\{\begin{array}{ll}i, & \mbox{if} \ i\in S(f) \\ i+c(f)-1, & \mbox{if}
\ i\not\in S(f). \end{array}\right.\]

The map $i\mapsto s_i$ is injective since $S(f)$ is a symmetric semigroup. On the
other hand, the map is surjective, because, given
$j\in S(f)\setminus (S(f)+c(f)-1)$, we have $j=s_j$ if  $j\leqslant c(f)-1$;
otherwise, if $j=i+c(f)-1$ for some $i>0$, then again by the symmetry of $S(f)$,
it follows that $j$ does not belong to $S(f)$ and therefore $j=s_i$.
\cqd

We believe that the converse of Theorem \ref{tame} is true, in the sense that if
$\mu(f)=c(f)$, then $S(f)$ is a tame semigroup, or, equivalently, if $p$ divides
any of the minimal generators of $S(f)$, then $\mu(f)>c(f)$. If this is so, we
would conclude from our result that if $\mu(f)=c(f)$, then $C_f$ is
$\mu$-stable.

To reinforce our conjecture, observe that the result of \cite{GB-P} proves it
when $\mult(f)<p$. The following example is a
situation where the converse holds and is not covered by the result in
\cite{GB-P}.

\begin{example} Let $p$ be any prime number and $n$ and $m$ two relatively prime
natural numbers such that $p\nmid n$. Then all curves given by
$f(X,Y)=Y^n-X^{mp}$ do not satisfy the condition $\mu(f)=c(f)$, since
$\mu(f)=\infty$ and $c(f)=(n-1)(mp-1)$. So, for all $p<n$, we have examples
for the converse of our result not covered by \cite{GB-P}.
\end{example}

Anyway, the other possible converse of Theorem \ref{tame}, namely, if $f$ is
$\mu$-stable then $S(f)$ is tame, is not true, as one may see from the following
example.

\begin{example} Let $f=(Y^2-X^3+X^2Y)^2-X^{11}Y\in k[[X,Y]]$, where $\car{k}=5$. Since $f^3\in \mathcal M T(f)^3$ (verified with Singular), then $C_f$ is $\mu$-stable, but its semigroup of values $S(f)=\langle 4,6,25\rangle$ is not tame.
\end{example}

\section{Proof of the Key Theorem}\label{ChapIII3}

This section is dedicated to prove Theorem 3.5, which proof is based on Proposition 4.2, that constitute the technical core of this work.

Let us start with an auxiliary result. Let $f\in k[[X,Y]]$ be an irreducible Weierstrass
polynomial in $Y$ of degree $n=v_0$, where
$S(f)=\langle v_0,\ldots, v_g\rangle$, $\I(f,X)= v_0$ and $\I(f,Y)= v_1$.

Consider the $k[[X]]$-submodule $V_{n-1}$ of $k[[X,Y]]$ generated by
$1,Y,\ldots,Y^{n-1}$,
and let $h_0=1, h_1, \ldots, h_{n-1}$ be polynomials in $Y$ such that
\[ V_{n-1} =  k[[X]]\oplus k[[X]]h_1 \oplus \cdots \oplus k[[X]]h_{n-1},\]
and their residual classes $y_i$ are the Ap\'ery generators of ${\mathcal O}_f$
as a free $k[[X]]$-module (cf. \cite[Proposition 6.18]{He}).

The natural numbers $a_i=v(y_i)=\I(f,h_i)$, $i=0,\ldots,n-1$, form the Ap\'ery
sequence of $S(f)$, so they are such that $0=a_0<a_1< \cdots <a_{n-1}$ and
$a_i \not\equiv a_j \bmod n$ for $i\neq j$ (cf. \cite[Proposition 6.21]{He}).

We have the following result.

\begin{proposition}\label{apery}
 Let $I$ be an $\mathcal M$-primary ideal of $k[[X,Y]]$ and $h\in V_{n-1}.$
If $\I(f,h)>>0$ then $h\in I$.
\end{proposition}

\pf Since the ideal $I$ is $\mathcal M$-primary, there exists a natural number
$l$ such that
$\mathcal M^l\subset I$.

Now, write $h=b_0+b_1h_1+\cdots+b_{n-1}h_{n-1}$,
with $b_i\in k[[X]]$, for all $i$. Since $\I(f,b_i) \equiv 0 \bmod n$,
$\I(f,h_i)=a_i$ and $a_i\not \equiv a_j \bmod n$, for $i,j=0, \ldots,n-1$, with
$i\neq j$, we have that
\[\I(f,h)=\min_i\{\I(f,b_i)+a_i\}\leqslant\min_i\{\I(f,b_i)\}+a_{n-1}.\]
Hence, $\I(f,h)>>0$ implies that $\min_j\{\I(f,b_j)\}>lv_0$. Therefore,
$h\in \mathcal M^l\subset I$, as we wanted to show.\cqd


Under the assumptions that $f$ is an irreducible Weierstrass polynomial in $Y$
with
$S(f)=\langle v_0,\ldots, v_g\rangle$ and $p\nmid v_0$, one may associate the
Abhyankar-Moh approximate roots  (cf. \cite[\S 6,7]{A-M}), which are
irreducible Weierstrass polynomials $f_{-1}=X,f_0=Y,\ldots,f_{g-1}$ such that,
for each $j\geqslant 0$, one has
$\deg_Y f_j=\frac{v_0}{e_j}$, $\I(f,f_j)=v_{j+1}$ and $S(f_j)= \langle \frac{v_0}{e_{j}}, \ldots, \frac{v_j}{e_{j}} \rangle$, satisfying a relation
\[f_j=f^{n_j}_{j-1}-\sum_{i=0}^{n_j-2}a_{ij}f_{j-1}^i,\]
where $a_{ij}$ are polynomials in $Y$ of degree less than
$\deg (f_{j-1})$ for $j=0,\ldots,g$ and
$\deg$ stands for the degree as a polynomial in $Y$.

So, from Corollary \ref{maisumdelgado} we have that
\begin{equation} \label{cordelgado}
\I(f,[f,f_{j-1}])\geqslant v_j+c(f)-1, \ \ \hbox{with equality if and only if} \
p\nmid v_j. \tag{4.1}
\end{equation}
This implies that if $p\nmid v_0v_1\cdots v_g$, then
$S(f)^*+c(f)-1\subset \nu(J(f)):=\{\I(f,h)\,\mid\,h\in J(f)\}$, where
$S(f)^*=S(f)\setminus \{0\}$.

The key result to prove Theorem \ref{theorem1} is Proposition
\ref{degreereduction} below that
will allow us to construct elements in $J(f)\cap V_{n-1}$ which intersection
multiplicity with $f$ sweep the set $S(f)^*+c(f)-1$.

\begin{proposition}\label{degreereduction} Let $f\in k[[X,Y]]$ be an irreducible
Weierstrass polynomial in $Y$ of degree $v_0$, where $k$ is an algebraically
closed field of characteristic
$p\geqslant 0$. Let $S(f)=\langle v_0,\ldots,v_g\rangle$ and suppose that
$p\nmid v_0v_1\cdots  v_g$. Given $s\in S(f)^*$, there exists $q_{s}\in J(f)$,
polynomial in $Y$, satisfying
\begin{enumerate}[{\rm (i)}]
\item $\deg {q_{s}}<\deg {f}=v_0$;
\item $\I(f,q_{s})=s+c(f)-1$.
\end{enumerate}
\end{proposition}
\pf We will use induction on $g$, the genus of $S(f)$, to construct step by step the
polynomial  $q_s$. It will be of the form $q_s=q_{f,s}=\sum_i P_i[f,f_{j_i}]$
(an infinite sum, possibly) where each $f_{j_i}$ is an approximate root of $f$
and the $P_i$ are monomials in the approximate roots of $f$ satisfying the
following conditions:
\begin{equation} \label{condition}
\left\{\begin{array}{llr} \I(f,P_1f_{j_1})=s; & \phantom{-}& \mbox{(4.2.1)}\\
\I(f,P_if_{j_i})>s, & \mbox{if} \ i\neq 1; & \mbox{(4.2.2)}\\
\deg {P_if_{j_i}}<\deg {f}, & \mbox{ for all } i; & \mbox{(4.2.3)}\\
\deg q_s < \deg f. & & \mbox{(4.2.4)}
\end{array} \right. \tag{4.2}
\end{equation}

This suffices to prove the proposition because (4.2.1) and (4.2.2) together with the
equality
$\I(f,P_i[f,f_{j_i}]) = \I(f,P_if_{j_i}) +c(f)-1$ (which follows from
(\ref{cordelgado})) imply (ii) in the statement.

If $g=0$, we have $f=Y$, so $J(f)=k[[X,Y]]$. Given
$s\in\mathbb{N}^*=S(f)^*$, set $$q_{f,s}:=X^{s-1}[f,X].$$
It is easy to check that $q_{f,s}$ satisfies (\ref{condition}).

Inductively, we assume that the construction was carried on for branches with genus $g-1$.
Consider the approximate root $f_{g-1}$ of $f$ which has genus $g-1$.
Since $n_gv_g\in \langle v_0,\ldots, v_{g-1} \rangle$ and $e_{g-1}=n_g$, we have
\[S(f)= \langle v_0,\dots,v_g\rangle \subset \Big\langle\frac{v_0}{n_{g}},\dots,
\frac{v_{g-1}}{n_{g}}\Big\rangle =S(f_{g-1}).\]

For $t\in S(f_{g-1})^*$, the inductive
hypothesis guarantees the existence of a $Y$-polynomial
\[q_{f_{g-1},t}=\sum_i P_i[f_{g-1},f_{j_i}],\]
where each $f_{j_i}$ is one of the approximate roots $f_{-1},f_0,\ldots,f_{g-2}$,
the $P_i$ are monomials in these approximate roots satisfying (\ref{condition}),
with $f_{g-1}$ and $v_0/e_{g-1}$ replacing $f$ and $v_0$, respectively.
Using this $q_{f_{g-1},t}$ we introduce the following
\textit{auxiliary polynomial}
\[\tilde{q}_{f_{g-1},t}:=\sum_i P_i[f,f_{j_i}].\]
To begin with, we will estimate the degree in $Y$ of these polynomials.
The inductive hypothesis gives $\deg {q_{f_{g-1},t}}<\deg {f_{g-1}}$ and
$\deg P_i\leqslant \deg {P_if_{j_i}}<\deg {f_{g-1}}$, for all $i$. On the other
hand, the Abhyankar-Moh's relation $f=f_{g-1}^{n_{g}}-G$, where
$G=a_{n_{g}-2}f_{g-1}^{n_{g}-2}+\dots+a_0$
and $\deg {a_{i}}<\deg {f_{g-1}}$, gives the inequality
$\deg {G} < (n_g-1)\deg {f_{g-1}} = \deg {f}-\deg {f_{g-1}}.$ We also have
$\deg {[G,f_{j_i}]}=\deg {(G_Xf_{j_i,Y}-G_Yf_{j_i,X}})\leqslant\deg {G}+
\deg {f_{j_i}}-1.$
It follows from the previous two inequalities and from the induction hypothesis
that
\[ \deg{P_i[G,f_{j_i}]} < \deg {P_i}f_{j_i}+ \deg {f}-\deg {f_{g-1}} -1
<\deg {f},\]
which together with the identity
\[\tilde{q}_{f_{g-1},t}=\sum_i P_i[f_{g-1}^{n_g}-G,f_{j_i}]=
n_gf_{g-1}^{n_g-1}q_{f_{g-1},t}-\sum_i P_i[G,f_{j_i}],\]
give the estimate
$$\deg{\tilde{q}_{f_{g-1},t}}<\deg{f}, \ \ \forall t\in S(f_{g-1})^*.$$\smallskip

\noindent\textbf{Claim:} \textit{For $t\in S(f_{g-1})^*$ we have
$\I(f,\tilde{q}_{f_{g-1},t})=c(f)-1+n_gt.$}\\
\\
Indeed, we have seen that $\I(f,P_i[f,f_{j_i}]) = \I(f,P_if_{j_i}) +c(f)-1.$
On the other hand, since the $P_if_{j_i}$ are products of approximate roots of
$f_{g-1}$ (so, also of $f$) and $\I(f_{g-1},f_i)=v_{i+1}/n_g$ we have
$\I(f,P_1f_{j_1})=\I(f_{g-1}^{n_g},P_1f_{j_1})=n_gt$.
Now, from (\ref{condition}), the intersection number $\I(f_{g-1},P_if_{j_i})$
assumes its minimum value once for $i=1$, when it is equal to $t$. Hence we have
$ \I(f,\tilde{q}_{f_{g-1},t}) = \I(f,\displaystyle\sum_i P_i[f,f_{j_i}])
= \I(f,P_1f_{j_1}) +c(f)-1 = n_gt+c(f)-1$, concluding the proof of the claim.
\medskip

The family of polynomials
$\{ \tilde{q}_{f_{g-1},t}; \ t \in S(f_{g-1})^* \}$
just introduced will be used in the construction of the family
$\{ q_{f,s}; \ \ s\in S(f)^*\}$ as announced in the beginning of the proof.

To this purpose, observe that each element $s$ of $S(f)^*$ decomposes
uniquely as
$$s=n_gt+wv_g, \ \ \hbox{with} \ t\in S(f_{g-1}), \ w\in \{0,1,\dots,n_g-1\}.$$

Now, we break up the analysis in three cases.\\
\\
{\bf Case $1$:} $(w=0)$ From Claim, we have
$s+c(f)-1=n_gt+c(f)-1=\I(f,\tilde{q}_{f_{g-1},t}).$
The estimate on the degree of $\tilde{q}_{f_{g-1},t}$, made just before Claim,
allows us to deduce that the series
\[q_{f,s}:=\tilde{q}_{f_{g-1},t}\]
has all the required properties, which proves (\ref{condition}) in this case.
\medskip

\noindent {\bf Case $2$:} $(t=0)$ In this case we prove conditions
(\ref{condition})
by induction on $w\geqslant 1$ with the extra and stronger condition
$\deg P_i f_{j_i} \leqslant w \deg f_{g-1}$ instead of $(2.3)$, since
$w\leqslant n_g -1$ and $\deg f_{g-1}=v_0/n_g$.

For $w=1$ the first two
conditions in
(\ref{condition}) and the extra condition hold trivially for
$$q_{f,v_g}:=[f,f_{g-1}].$$
Moreover, since $f=f_{g-1}^{n_g} -\, G$ by the preceding estimates we get
$\deg{q_{f,v_g}}  =  \deg {[f_{g-1},G]}
  \leqslant \deg{G}+\deg{f_{g-1}}-1
  < (\deg{f}-\deg{f_{g-1}})+\deg{f_{g-1}}-1
  <  \deg{f}.$

For the induction step we need the following result that gives a method to
reduce degrees while preserving
intersection multiplicities with $f$ and residual classes modulo
$J(f)$ and which proof will be given later in order to not interrupt the proof
of Proposition 4.2 in course.

\begin{lemma} \label{reducaodograu} Let $m\in\mathbb{N}^*$ be such
that $n_g\nmid m$ and $m>n_g(c(f_{g-1})-1)$. Suppose that we have a
$Y$-polynomial $h$ such that
$\deg{h}<\deg{f} \ \ \hbox{and} \ \ \I(f,h)=c(f)-1+m.$ Then there  exists a
$Y$-polynomial $h'$, such that
\begin{itemize}
 \item[(i)] $\I(f,h')=\I(f,h)$;
 \item[(ii)] $\deg{h'} < \deg{f}-\deg{f_{g-1}}$;
 \item[(iii)] $h-h'=\displaystyle\sum_{j}\alpha_j\tilde{q}_{f_{g-1},u_j} ,
 \ \  \alpha_j\in k, \ \ u_j\in  S(f_{g-1}), \ \ n_gu_j>m, \forall j.$
\end{itemize}

\end{lemma}

To continue our proof of Proposition 4.2, we are going to apply the above Lemma to $m=(w-1)v_g$.

Notice that
$(w-1) v_g \geqslant v_g>n_g(c(f_{g-1})-1)$. Indeed to see the last inequality,
since $S(f)$ is strongly increasing we have $v_i-v_{i-1}>(n_{i-1}-1)v_{i-1}$.
Summing up all these inequalities, for $i=1,\ldots, g$, with $n_0:=1$ we obtain
$v_g-v_0>c(f)-1 - (n_g - 1)v_g = n_g(c(f_{g-1})-1)$, were the equality comes from
conductor formula.
On the other hand, since gcd$(v_g,n_g)=1$, we have that $n_g\nmid m$.

Let us suppose that $q_{f,(w-1)v_g}$ is already constructed satisfying
(\ref{condition}) and the extra condition.

By using the previous lemma for $h=q_{f,(w-1)v_g}$, we define
\[ q_{f,wv_g} = f_{g-1} \left( q_{f,(w-1)v_g} \right)'.  \]

In order to check (4.2.1) and (4.2.2) it suffices to notice that in the expansion of
$\left( q_{f,(w-1)v_g} \right)'$ there exists only one index $i=1$ such that
$\I(f, P'_1f_{j_1})=(w-1)v_g$ and $\I(f, P'_if_{j_i})>(w-1)v_g$ for all $i\neq 1$.

The extra condition on the degrees is obtained by observing that for the expansion of
$q_{f,(w-1)v_g}$ we have $\deg P_if_{j_i} \leqslant (w-1)\deg f_{g-1}$ and for
the expansion of
$\tilde{q}_{f_{g-1},u_j}$ we have
$\deg \tilde{P}_if_{j_i} < \deg f_{g-1} \leqslant (w-1)\deg f_{g-1}$.

The fourth condition (4.2.4) follows from (ii) in the above Lemma. This concludes
the proof of this case.

\medskip

\noindent {\bf Case $3$:} $(w>0$ and $t>0)$
From Case 2, we obtain $q_{f,wv_g}$ satisfying conditions (\ref{condition})
with the extra condition on degrees. From the induction hypothesis on the genus
we get $q_{f_{g-1},t}=\sum_i P''_i[f_{g-1},f_{m_i}]$ satisfying (\ref{condition})
for $f_{g-1}$ instead of $f$. Now we define
\[ q_{f,s}= P''_1f_{m_1}\left( q_{f,wv_g} \right)'=
P''_1f_{m_1}\left(q_{f,wv_g}+
\displaystyle\sum_{j}\alpha_j\tilde{q}_{f_{g-1},u_j}\right).\]
In order to check (4.2.1) and (4.2.2) it suffices to notice that
$\I(f,P''_1f_{m_1})=n_gt$, as we have seen in the proof of Claim, and
that in the expansion of
$\big( q_{f,wv_g} \big)'$ there exists only one index $i=1$ such that
$\I(f, P'_1f_{j_1})=wv_g$ and $\I(f, P'_if_{j_i})>wv_g$ for all $i\neq 1$.

In order to prove (4.2.3) remind that $\deg P''_1f_{m_1} < \deg f_{g-1}$, that
for the expansion of
$q_{f,wv_g}$ we have $\deg P_if_{j_i} \leqslant w\deg f_{g-1}$ and for
the expansion of
$\tilde{q}_{f_{g-1},u_j}$ we have
$\deg \tilde{P}_if_{j_i} < \deg f_{g-1} \leqslant w\deg f_{g-1}$. Therefore,
for terms $P'_if_{j_i}$ in the expansion of $\big( q_{f,wv_g} \big)'$ we have
$\deg P''_1f_{m_1}P'_if_{j_i} < \deg f_{g-1} + w\deg f_{g-1}=(w+1)\deg f_{g-1}
\leqslant (n_g-1+1)\deg f_{g-1}=\deg f$.

Condition (4.2.4) follows directly from (ii) in Lemma \ref{reducaodograu} and from
$\deg P''_1f_{m_1} < \deg f_{g-1}$.

To finish the proof of the Proposition \ref{degreereduction} it remains to prove
Lemma \ref{reducaodograu}.




\noindent{\bf Proof of Lemma \ref{reducaodograu}:}
Indeed, we begin by dividing $h$ by $f_{g-1}^{n_g-1}$. Then we get
$h=f_{g-1}^{n_g-1}h_0''+h_0'$ where $\deg{h_0'}<\deg{f_{g-1}^{n_g-1}}=
\deg{f}-\deg{f_{g-1}}.$ The rough idea of the proof is to eliminate the term
$f_{g-1}^{n_g-1}h_0''$ in the preceding relation using the polynomials
$\tilde{q}_{f_{g-1},u}$ where $u\in S(f_{g-1})^*$. This will be done iteratively,
in possibly infinitely many steps, with the help of the following auxiliary
result, which we prove after we finish the proof of Lemma \ref{reducaodograu}.

\begin{lemma} \label{auxiliaryresult}
With the same conditions as above, we have
$\I(f,h_0'')=n_g\I(f_{g-1},h_0'')$  and
$\I(f,h_0''f_{g-1}^{n_g-1})\neq \I(f,h_0').$
\end{lemma}

Using the formula $c(f)-1=n_g(c(f_{g-1})-1)+(n_g-1)v_g$ and Lemma
\ref{auxiliaryresult}, we get
$\I(f,h_0''f_{g-1}^{n_g-1})-(c(f)-1)=n_g[\I(f_{g-1},h_0'')-c(f_{g-1})+1].$
On the other hand, since \linebreak
$\I(f,h)-(c(f)-1)=m$ and $n_g\nmid m$, it follows that
 \[\I(f,h_0')=\I(f,h)<\I(f,f_{g-1}^{n_g-1}h_0'').\]
So, from the first part of Lemma
\ref{auxiliaryresult} and the above inequality, we get
$$n_g\I(f_{g-1},h_0'')=\I(f,h_0'')>\I(f,h)-\I(f,f_{g-1}^{n_g-1})=
m+c(f)-1-(n_g-1)v_g.$$

Defining $u_1=\I(f_{g-1},h_0'')-c(f_{g-1})+1$, it follows that
$c(f_{g-1})-1+u_1= \I(f_{g-1},h_0'')>\frac{m}{n_g}+c(f_{g-1})-1>2(c(f_{g-1})-1),$
allowing us to conclude that $u_1\in S(f_{g-1})^*$.

The inductive hypothesis guarantees the existence
of a polynomial $q_{f_{g-1},u_1}$ satisfying all requirements in
(\ref{condition}).

From Claim, we have
$\I(f,\tilde{q}_{f_{g-1},u_1})=c(f)-1+n_gu_1=\I(f,h_0''f_{g-1}^{n_g-1}).$
So, after multiplication by a suitable $\alpha_1\in k^{*}$, we get that
$h_1:=h_0''f_{g-1}^{n_g-1}-\alpha_1\tilde{q}_{f_{g-1},u_1}$ satisfies
$\deg h_1 < \deg f$ and the
inequality
\begin{equation} \label{firststep} \I(f,h_1)> \I(f,h_0''f_{g-1}^{n_g-1})>\I(f,h) \
=c(f)-1+m.
\tag{4.3}
\end{equation}

This allows us to write
\[
h=h_1+\alpha_1 \tilde{q}_{f_{g-1},u_1}+h_0', \ \ \hbox{with} \ \
\I(f,h_1)>\I(f,h) \ \  \hbox{and} \ \ \I(f,h_0')=\I(f,h).
\]
From (\ref{firststep}) and from the definition of $h_1$ we have that there exists
$m_1\in {\mathbb N}^*$ such that
\[
\I(f,h_1)=c(f)-1+m_1 >c(f)-1 +n_gu_1>c(f)-1+m.
\]
So, $m_1>m$ and $n_gu_1>m$.

In the next step we proceed differently according to the divisibility of $m_1$
by $n_g$.

Suppose $n_g\nmid m_1$. Hence we are in position to repeat the preceding
process of division by $f_{g-1}^{n_g-1}$ using, this time, $h_1$ instead of $h$.
So $h_1=f_{g-1}^{n_g-1}h_1''+h_1'$. Again,  we deduce that there exist
$\alpha_2\in k^*$ and $u_2\in S(f_{g-1})^*$, with $n_gu_2>m_1>m$, such that if
we define $h_2:=h_1''f_{g-1}^{n_g-1}-\alpha_2\tilde{q}_{f_{g-1},u_2}$, then we
have $\deg h_2 < \deg f$ and
$$\I(f,h_2)>\I(f,h_1)>\I(f,h).$$

If, however, $n_g\mid m_1$, say $m_1=n_gu_2$, by the inequality just after
the definition of $u_1$ and the inequality that follows from (\ref{firststep})
we have
\[c(f_{g-1})-1+u_2>c(f_{g-1})-1+u_1>2(c(f_{g-1})-1).\] So, it follows that
$u_2\in S(f_{g-1})^*$. Hence, there exists a polynomial
$q_{f_{g-1},u_2}$ such that
\[\I(f,h_1)=\I(f,\tilde{q}_{f_{g-1},u_2} ) \] and again we may choose
$\alpha_2\in k^*$ in such a way that if
$h_2:=h_1-\alpha_2\tilde{q}_{f_{g-1},u_2}$, we have $\deg h_2 < \deg f$ and
$\I(f,h_2)>\I(f,h_1)$.
Hence, we get
$h=h_2+\alpha_1\tilde{q}_{f_{g-1},u_1}+h_0'+\alpha_2\tilde{q}_{f_{g-1},u_2}+h_1'$,
where $h_1'=0$, in this case. Notice that $n_gu_2>n_gu_1>m$.

So, by repeating this process we obtain
$$h=h_j+ \sum_{i=1}^j \alpha_i\tilde{q}_{f_{g-1},u_i}+ \sum_{i=0}^{j-1}h_i',$$
with $\I(f,h_i')<\I(f,h_{i+1}')$ if $h_i'\neq 0$,
$\I(f,\alpha_i\tilde{q}_{f_{g-1},u_i})<
\I(f,\alpha_{i+1}\tilde{q}_{f_{g-1},u_{i+1}})$
and $\I(f,h_i)<\I(f,h_{i+1})$.
Since all power series appearing in the above sum have degree less than $\deg{f}$,
it follows, in view of Proposition \ref{apery}, that $h_j\rightarrow 0$ in the
$\mathcal M$-adic topology of $k[[X,Y]]$ and the families
$\{h_i'\}_{i\in {\mathbb N}}$ and
$\{ \alpha_i\tilde{q}_{f_{g-1},u_i} \}_{i\in {\mathbb N}} $ are
summable. Taking $h'=\sum_j h_j'$ we get Lemma \ref{reducaodograu}.\\




Finally, it remains to prove Lemma \ref{auxiliaryresult} in order to finish the
proof of Lemma \ref{reducaodograu} and finally conclude the proof of
Proposition \ref{degreereduction}.\\

\noindent{\bf Proof of Lemma \ref{auxiliaryresult}:}
If $f$ is any irreducible Weierstrass polynomial of degree $n$, then it is easy
to see from Proposition \ref{apery} that the set $V_{n-1}$ of all polynomials in
$Y$ of degree less than $n$ with coefficients in $k[[X]]$ is a free
$k[[X]]$-module with basis
$$\left\{f^J=f_0^{j_0}f_1^{j_1}\cdots f_{g-1}^{j_{g-1}}; \ J=(j_0,\ldots,j_{g-1}),
\ 0\leqslant j_i<n_{i+1}, \ i=0,\ldots,g-1 \right\}.$$ So, any element
$h\in V_{n-1}$ is written uniquely as
$h=\sum_J a_J(X)f^J =
f_{g-1}^{n_g-1}h''+h'$, where $a_J(X)\in k[[X]]$ and

\[
h''=\displaystyle\sum_{j_{g-1}=n_g-1}a_J(X)f_0^{j_0}
\cdots f_{g-2}^{j_{g-2}} \ \ \ \mbox{ and } \ \ \ h'=\displaystyle\sum_{j_{g-1}\leqslant
n_{g}-2}a_J(X)f_0^{j_0}\cdots f_{g-1}^{j_{g-1}}.
\]

First of all we will check that $\I(f,h')\neq \I(f,f_{g-1}^{n_g-1}h'')$. In fact,
it follows, from the uniqueness of the canonical representation of the elements of $S(f)=\langle v_0,\ldots,v_g\rangle$,
that in $h'$ there is a unique term such that
\[\I(f,h')=\I(f,a_J(X)f_0^{j_0}\cdots f_{g-1}^{j_{g-1}})=
\sum_{i=-1}^{g-1}j_iv_{i+1},\]
where $j_{-1}=\ord_Xa_J(X)$ and $j_{g-1}\leqslant n_g-2$. Also, in
$f_{g-1}^{n_g-1}h''$ there is a unique term satisfying
\[\I(f,f_{g-1}^{n_g-1}h'')=
\I(f,a_L(X)f_0^{l_0}\cdots f_{g-2}^{l_{g-2}}f_{g-1}^{n_g-1})=
\sum_{i=-1}^{g-2}l_iv_{i+1}+(n_g-1)v_g, \]
where $ l_{-1}=\ord_X(a_L(X))$. Now the claim stated above follows from the
uniqueness of writing in $S(f)$.

Also, it is clear from the way we wrote $h''$ that
$\I(f,h'')=n_{g}\I(f_{g-1},h'')\in n_gS(f_{g-1})$.

Now, to finish the proof of this Lemma we only need to check that $h''$
and $h'$ are indeed the quotient and the remainder, respectively, of the
division of
$h$ by $f_{g-1}^{n_g-1}$. We will do this by estimating the degree of $h'$ and,
hence,
conclude by the uniqueness of the remainder and the quotient in the euclidean
algorithm.
Indeed, for every summand in $h'$ we have
$\deg{(a_J(X)f_0^{j_0}\cdots f_{g-1}^{j_{g-1}})}< \deg{f_{g-1}^{n_g-1}}$, which
shows that  $\deg{h'}<\deg{f_{g-1}^{n_g-1}}$.\cqd \smallskip



With these tools at hands, we may conclude the proof of Theorem \ref{theorem1}.\\

\noindent \textit{\bf Proof of the Theorem \ref{theorem1}:} Choose ${\mathcal{F}}$ with minimal number of elements, so from Lemma \ref{contagem} it follows that $\#{\mathcal{F}}=c(f)$. We will show that
the set $\overline{\mathcal{F}}$ generates $k[[X,Y]]/J(f)$ as a $k$-vector space.
In particular, this will show also that $\mu(f)\leqslant c(f)$ when $S(f)$ is tame.
In order to do this, it is enough to show that there exists a decomposition
$k[[X,Y]]=\langle \mathcal{F}\rangle + J(f),$ where $\langle \mathcal{F}\rangle$ denotes the $k$-vector space spanned by the elements of $\mathcal{F}=\{\varphi_1,\ldots,\varphi_{c(f)}\}$.

Given any element $h\in k[[X,Y]]$ we can divide it  by the partial derivative $f_Y$
which, under our assumptions, is a $Y$-polynomial of degree $v_0-1$.
The remainder of the division is a $Y$-polynomial $h'$ of degree less than $v_0-1$
and it is sufficient to show that $h'$ belongs to
$\langle \mathcal{F}\rangle + J(f).$

If $\I(f,h')\in S(f)\setminus (S(f)+c(f)-1)$ then, according to the definition of $\mathcal{F}$, there is an element $\varphi_{i_{s_0}}$ such that
$s_0:=\I(f,h')=\I(f,\varphi_{i_{s_0}})$. Hence, there is a constant $\alpha_{s_0}\in k$ such that
\[\I(f,h'-\alpha_{s_0}\varphi_{i_{s_0}})=:s_1>s_0.\]

If otherwise $\I(f,h')=s_0=s'_0+c(f)-1\in S(f)^*+c(f)-1$, then Proposition \ref{degreereduction} guarantees the
existence of an element $q_{s'_0}$ in $J(f)$, polynomial in $Y$ of degree less than $\deg{f}$, such that $s_0=\I(f,q_{s'_0})$. Hence,
there is a constant $\beta_{s_0}\in k$ such that
\[\I(f,h'-\beta_{s_0}q_{s'_0})=:s_1>s_0.\]

We carry on this process that increases intersection indices to eventually achieve
\[s_r=\I\left(f,h'- \sum_s\beta_{s}q_{s'}-  \sum_s \alpha_s\varphi_{i_s}\right)\in S(f)^*+c(f)-1, \ \ \ \mbox{for all}\,\, r\geqslant N.\] Since the elements in $S(f)^*+c(f)-1$ may be realized as intersections indices of $f$ with elements in $J(f)\cap V_{n-1}$ (cf. Proposition \ref{degreereduction}), we produce an element \[h'-\sum_s \alpha_s\varphi_{i_s}-\sum_s \beta_s q_{s'}\]
which intersection multiplicity with $f$ is big enough and which degree is less than $\deg{f}$. We aim to use
Proposition \ref{apery} to conclude that it belongs to the Jacobian ideal $J(f)$. In order to do so we need only to check that
$J(f)$ is a $\mathcal M$-primary ideal, as we did in Lemma \ref{tameimplicafinito}.
\cqd


As a final comment, we recall the classical Milnor's Formula for plane curves
singularities that states that if $f=f_1\cdots f_r\in \C[[X,Y]]$ is a possibly
many branched and reduced power series over $\C$ then $\mu(f)=2\delta(f)+1-r$.
A first and natural question after obtaining Theorem \ref{tame} is:

\medskip

\noindent \emph{If all branches $f_1,\ldots,f_r$ have tame semigroups does the
preceding formula continues to hold?}

\medskip

We include here some examples to show that this is not always true.

\medskip

\begin{example}
Let $f=(Y^2-X^3)^2-X^{11}Y$ and $g=(Y^2-X^3+X^2Y)^2-X^{11}Y$. Then $\I(f,g)=28$
and $S(f)=S(g)=\langle4,6,25\rangle$. Here one can compute
$\delta(fg)=\delta(f)+\delta(g)+\I(f,g)=14+14+28=56.$
Hence $2\delta(fg)+1-r=111$. The behaviour of the reduced $2$-branched plane
curve singularity defined by the equation $fg$ with respect to $\mu$-stability
and the value of
$e_0(T(fg))=\mu(\Oh_{fg})$ according to the characteristic $p$ of the ground
field is described below. All computations are performed using the software
\emph{Singular}, \cite{DGPS}.\smallskip

\noindent{\bf $(p=7)$} Here $(fg)^4\in\M\,T(fg)^4$ so that $fg$ is $\mu$-stable. Computation shows that $\mu(\Oh_{fg})=\mu(fg)=112>111$ so that Milnor's Formula does not hold. Notice that the semigroups of the branches are tame in characteristic $p=7$, but $p\mid \I(f,g)$.

\noindent {\bf $(p=5)$} Here $(fg)^3\in\M\,T(fg)^3$ so that again $fg$ is $\mu$-stable. Computation now shows that $\mu(\Oh_{fg})=\mu(fg)=111$ and Milnor's Formula does hold. Notice that the semigroups of the branches are wild in characteristic $p=5$.

\noindent {\bf $(p=13)$} Here computation shows that $\mu(fg)=124$ and
 $\mu((1+X)fg)=114$ so that $fg$ is not $\mu$-stable. Computation also suggests that $\mu(\Oh_{fg})=\mu(fg)=114$ and Milnor's Formula does not hold. Notice that the semigroups of the branches are tame in characteristic $p=13$ and $p\nmid \I(f,g)$.
\end{example}




\noindent Authors addresses: {ahefez@id.uff.br; \ joaohelder@id.uff.br; \ rsalomao@id.uff.br}




\end{document}